\RequirePackage{fix-cm}
\documentclass{svjour3}                     
\smartqed  

\usepackage{graphicx,amsmath,amssymb,enumerate,mathrsfs,calrsfs,lineno,hyperref}
\usepackage{caption}
\usepackage{subcaption}

\begin{document}

\title{ A Fractional Image Inpainting Model Using a Variant of Mumford-Shah Model  
}

\author{Abdul Halim       \and
         Abdur Rohim \and
         B.V. Rathish Kumar \and
         Ripan Saha
}


\institute{\\ Department of Applied Mathematics \\
King Abdullah University of Science and Technology, Saudi Arabia\at
              \email{abdul.math91@gmail.com}           
           \and
           {Department of Mathematics, Raiganj University, West Bengal, India}\\
        \email{abdurmath90@gmail.com}    \\
         \and
           {Department of Mathematics and Statistics , IIT Kanpur, Uttar Pradesh, India}\\
        \email{drbvrk11@gmail.com}  
        \and
           {Department of Mathematics, Raganj University, India}\\
}

\date{\today}

\maketitle

\begin{abstract}
In this paper, we propose a fourth order PDE model for image inpainting based on a variant of the famous Mumford-Shah (MS) image segmentation model. Convexity splitting is used to discrtised the time and we replace the Laplacian by its fractional counterpart in the time discretised scheme. Fourier spectral method is used for space discretization. Consistency, stability and convergence of the time discretised model has been carried out. The model is tested on some standard test images and compared them with the result of some models existed in the literature.

\keywords{Image Processing \and Inpainting \and Fourth Order Partial differential equation (PDE)\and Convexity splitting.}
\end{abstract}

\section{Introduction}
\label{intro}
Image inpainting is a process of filling in the missing part of an image from the information available outside the missing part. In mathematical term it  is an extrapolation of the image. Although there are several inpainting methods are available in the literature, we consider the Partial differential equation (PDE) based methods because of its sound established theory.  \\

Since image inpainting is relatively new topic in the field of image processing, a lot of inpainting models are based on the models of denoising and segmentation. Like in  \cite{CS1} Chan and Shen have used a modified version of Rudin-Osher-Fatemi (ROF) denoising model \cite{rof} for image inpainting and the resulting inpainting model is known as $TV-L^2$ inpainting model. Following the $TV-L^2$  model a series of models \cite{sch,CS2,CS3} have been proposed which fall in the category of variational inpainting. In variational inpainting, the PDE is obtained by minimizing an energy functional of the following form:
 	\begin{equation}\label{VF}
 	E(u)=R(u)+\lambda F(f-u)
 	\end{equation}
      where $f$ is the given damaged image and $u$ is the inpainted image and
      \begin{equation}\label{fid}
       \lambda (x)=\begin{cases}
                    \lambda_0 \hspace{.1cm} in \hspace{.2cm}\Omega \setminus D\\
                    0 \hspace{.2cm}in \hspace{.2cm} D
                   \end{cases}
      \end{equation}
      with $\lambda_0 >> 1$ and $D$ is the damaged part of the image.
The term $R(u)$ is called regularization term, $F(\cdot)$ is called fidelity term and $\lambda$ is called fidelity parameter which forces the inpainting image to remain closer to the original image outside the missing/ degraded region $D$.
    
Another variational model for inpainting has been proposed by Esedoglu and Shen \cite{mse}. They have brought the Mumford-Shah model of segmentation and used it for image inpainting. The Mumford-Shah(MS) functional is non-convex and contains an unknown set of lower dimension, so it is difficult to solve the minimization problem. Several approximations have been proposed in the literature~\cite{surveyMS}. Esedoglu and Shen \cite{mse} have used the Ambrosio and Tortorelli’s approximation of the MS energy functional for the inpainting model. Applying the gradient descent to the approximations give raise to a parabolic problem of second order with a parameter $\epsilon$. Since the second order models are not good for inpainting~\cite{tvh}, as second order models are unable to fill large gaps, they have modified their second order model using Euler-Elastica energy to get a fourth order model in the same paper \cite{mse}, known as Mumford-Shah-Euler (MSE) model.

The parameter $\epsilon$ of the MSE model reminds the authors of \cite{mCH} about the Cahn-Hilliard equation used to model phase-separation phenomena of material science. They have modified the Cahn-Hilliard equation and proposed a model for image inpainting in \cite{mCH}. Here after we will call this model as mCH model. The mCH model has few drawbacks, like it is applicable for binary images only. Then several authors have come up with inapinting models \cite{cCH,fCH,anis,lininp,nonlininp} by modifying the Cahn-Hilliard model to overcome the drawbacks of mCH model.

More recently, Cai et.al.~\cite{twophase} have proposed a segmentation model which is also derived from the Mumford-Shah model. This segmentation model is also used for inpainting in the subsequent paper of the same authors in ~\cite{threephase}. But this model will leads to second order equation. As we know that the second order model are not able to fill large gaps~\cite{tvh}, so we will propose a fourth order model in this article using a variant of the MS functional. We will use convexity splitting~\cite{tvh} in time and Fourier spectral methods in space for the proposed model.

The paper is organised is as follows: In section \ref{sec:models} we have discussed some relevant existing models and thereby proposed our model. Section \ref{numschm}, talks about convexity splitting and the numerical scheme for our proposed model. Also we have established the consistency, stability and convergence of the time discretised model. Further, we have introduced the fractional time discretised model. In section \ref{numresult}, we have presented the fractional version of the time discretized scheme and present the scheme with complete discretization. Then the numerical results of our model have been presented and  compared them with the same of models like $TV-L^2$ and $TV-H^{-1}$ and the second order model of Cai et.al. \cite{twophase}. Finally, in section  \ref{conclusion} we draw some conclusions.
 	       
\section{Proposed model and some existing models}
\label{sec:models}
In this section, we will discuss few inpainting models from the literature which are relevant to us and then proposed our model.
 \subsection{$TV-L^2$ and $TV-H^{-1}$ model}
 The inpainting model $TV-L^2$ and $TV-H^{-1}$ are of variational type that is they are obtained by minimizing an energy functionl of the form \eqref{VF}. For both the models the regularization term is the total variation of $u$ that is $ R(u)=\int_{\Omega}| \nabla u | dx$. The fidelity term for the $TV-L^2$ model is $\|f-u\|_2$ and  $TV-H^{-1}$ is $\|f-u\|_{H^{-1}(\Omega)}$.
 
 With the above choices, the minimizing energy of the $TV-L^2$ model is :
 \begin{equation}\label{tvl2}
 	 E(u)=\int_{\Omega}|\nabla u|dx +\frac{\lambda}{2}\int_{\Omega}(f-u)^2dx.
 	\end{equation}
The gradient descent of the corresponding Euler-Lagrange equation is
 \begin{equation}\label{tvl2feq1}
         u_t=\nabla\cdot\frac{\nabla u}{|\nabla u|}+\lambda (f-u).
        \end{equation}
 To avoid division by zero the above equation is modified~\cite{CS1} as:
 \begin{equation}\label{tvl2feq}
         u_t=\nabla\cdot\frac{\nabla u}{\sqrt{|\nabla u|^2+\delta^2}}+\lambda (f-u),
        \end{equation}
        where $\delta<<1$ is a parameter.
       
Similarly, we get the PDE for the $TV-H^{-1}$ model as
        \begin{equation}\label{tvhfeq}
         u_t=-\Delta \Big( \nabla\cdot\frac{\nabla u}{\sqrt{|\nabla u|^2+\delta^2}}\Big)+\lambda (f-u),
        \end{equation}
        
   \subsection{Mumford-Shah model}
 As we have mentioned in the introduction that our model will be derived from MS-model. In MS- model on has to minimize the following energy functional
 \begin{equation}\label{eMS}
 E_{MS}(u,\Gamma)=\frac{\eta}{2}\int_{\Omega}(f-u)^2 dx+\frac{\mu}{2}\int_{\Omega\setminus \Gamma}|\nabla u|^2 dx +Length(\Gamma),
 \end{equation}
 where $\eta,\mu$ are parameters and $\Gamma$ is a curve inside the image domain $\Omega.$
 But the functional $E_{MS}$ is non-convex, so it is not solvable directly. So there are several simplified version has been proposed like proposed in \cite{twophase,mse}. A detail survey on Mumford-Shah model and its variant applied to image processing can be found in \cite{surveyMS}. 
 
 In \cite {twophase} authors have proved that minimizing the functional \eqref{eMS} is equivallent to minimizing the following functional:
  \begin{equation}\label{etMS}
 E(u)=\frac{\eta}{2}\int_{\Omega}(f-u)^2 dx+\frac{\mu}{2}\int_{\Omega}|\nabla u|^2 dx + \int_{\Omega}|\nabla u | dx
 \end{equation}
 For segmentation $\eta$ is constant throughout the domain but for inpainting $\eta$ is some constant outside the damaged part $D$ and zero in $D$. Changing the fidelity parameter $\eta$ of the model \eqref{etMS}, the same model has been used in \cite{threephase} in context of colour image which is of the form:
  \begin{equation}\label{etMS1}
 E(u)=\frac{\lambda}{2}\int_{\Omega}(f-u)^2 dx+\frac{\mu}{2}\int_{\Omega}|\nabla u|^2 dx + \int_{\Omega}|\nabla u | dx
 \end{equation}
 where $\lambda$ is defined in \eqref{fid}.
 The corresponding Euler-Lagrange equation will be
 \begin{equation}
 -\lambda(f-u)-\mu \Delta u -\nabla \cdot \frac {\nabla u}{|\nabla u|}=0.
 \end{equation}
 which is equivalent to solve the steady state solution of the equation 
  \begin{equation} \label{vms}
u_t= \mu \Delta u +\nabla \cdot \frac {\nabla u}{|\nabla u|}+\lambda(f-u).
 \end{equation}
 Here after we will call this model as convex variant of Mumford-Shah (CVMS) model.
  \subsection{Proposed Model}
 But this is a second order model so it will not be able to fill large gap. So we will propose a higher order model. To do this we modify the MS-functional as follows:
\begin{equation}\label{meMS}
 E_{MS}(u,\Gamma)=\frac{\eta}{2}\int_{\Omega}(f-u)^2 dx+\frac{\mu}{2}\int_{\Omega\setminus \Gamma} (u_{x_1x_1}^2+2u_{x_1x_2}^2+u_{x_2x_2}^2)dx +Length(\Gamma),
 \end{equation}
 then following the steps of \cite{twophase} we can prove that this is equivalent to minimising
   \begin{equation}\label{emMS}
 E(u)=\frac{\lambda}{2}\int_{\Omega}(f-u)^2 dx+\frac{\mu}{2}\int_{\Omega}|\Delta u|^2 dx + \int_{\Omega}|\nabla u | dx.
 \end{equation}
 Now using the regularized version of the total-variation the functional reduced to:
 \begin{equation}\label{emMS1}
 E(u)=\frac{\lambda}{2}\int_{\Omega}(f-u)^2 dx+\frac{\mu}{2}\int_{\Omega}|\Delta u|^2 dx + \int_{\Omega}{\sqrt{|\nabla u|^2+\delta^2}}.
 \end{equation}
 where $\delta$ is a small parameter. For all numerical results presented in this work we choose $\delta = 0.01$.   
The corresponding descent model will be 
  \begin{equation}\label{prop2}
 u_t=-\mu \Delta^2 u +\nabla \cdot \frac {\nabla u}{\sqrt{|\nabla u|^2+\delta^2}} +\lambda(f-u).
 \end{equation}
 
    \begin{table}[t]
 	\centering
 	\begin{tabular}{|c|c|c|c|c|c|c|c|} 
 		\hline
 		 Image & {\begin{tabular}[c]{@{}c@{}} $\alpha$\end{tabular}} &
 		 {\begin{tabular}[c]{@{}c@{}} PSNR \end{tabular}} & 
 		{\begin{tabular}[c]{@{}c@{}} SNR \end{tabular}}& 
		{\begin{tabular}[c]{@{}c@{}} SSIM \end{tabular}} & 
 		{\begin{tabular}[c]{@{}c@{}} Iteration \end{tabular}}
 		&   {\begin{tabular}[c]{@{}c@{}} Time taken \end{tabular}}\\	
\hline
Dog & 1 & 37.50&	30.46& 0.9412 &	2999&  17.29 \\
   &1.2& 38.78&	31.73&	0.9704&	1580&9.13 \\
   &1.4& 39.14&	32.09&	0.9767&	950 &5.41\\
   & 1.6& 39.13&	32.09&	0.9764&	830 &4.80\\
   &1.8&39.02 &	31.97&	0.9758&	810&4.74 \\
   &2& 38.88&	31.83&	0.9748&	720& 4.21  \\
\hline
GrayShade & 1 &35.11&	31.41&	0.7853	&1180& 9.10  \\
   &1.2& 42.39	&38.69&	0.9627&	1200& 9.37\\
   &1.4& 44.73	&41.03&	0.9872&	1210& 9.16\\
   & 1.6& 45.39&	41.69&	0.9925&	1220& 9.10\\
   &1.8&45.40	&41.70&	0.9925&	1170&9.15\\
   &2& 45.38&41.67&	0.9925&1200 &9.42\\
\hline
 \end{tabular}
 \caption{Comparison of results with different fractional power $\alpha$, in terms of PSNR, SNR and SSIM.}
 	\label{table3}
 \end{table}
 
\begin{table*}[t]
 	 	\centering
 	\begin{tabular}{|c|c|c|c|c|c|c|c|} 
 		\hline
 		 Image & {\begin{tabular}[c]{@{}c@{}} Model \end{tabular}} &
 		 {\begin{tabular}[c]{@{}c@{}} PSNR \end{tabular}} & 
 		{\begin{tabular}[c]{@{}c@{}} SNR \end{tabular}}& 
 		{\begin{tabular}[c]{@{}c@{}} SSIM \end{tabular}} & 
 		{\begin{tabular}[c]{@{}c@{}} Iteration \end{tabular}}
 		&   {\begin{tabular}[c]{@{}c@{}} parameters \end{tabular}}
		&   {\begin{tabular}[c]{@{}c@{}} CPU time \end{tabular}}\\

      \hline
   GrayShade &$TV-L^2$  &41.64 &	37.94&	0.9822&	33965&  $\lambda$=10& 18.81\\
                       &$TV-H^{-1}$& 45.23&	41.54&	0.9909 &	6576 & $\lambda$=10& 38.84 \\
                      &Our& 45.40&	41.70&	0.9925&	1210 &$\alpha$=1.8&  8.90\\
                  \hline
                  
                Kaleidoscope &$TV-L^2$  & 26.24&	18.51&	0.8827&	66365 &  $\lambda$=50& 52.35\\
                       &$TV-H^{-1}$& 26.90&	19.17&	0.8770&	2280& $\lambda$=50 & 22.87 \\
                        &Our& 27.89& 	20.16&	0.9130&	2020& $\alpha$=1.8& 23.07\\
                          \hline 
       Elephant 1 &$TV-L^2$  &35.02	&28.20	&0.9694&	2720&	 $\lambda$=200 &12.47 \\
   &$TV-H^{-1}$&33.93&	27.12&	0.9396	&2030&  $\lambda$=200 &11.41     \\
      &Our&36.16&	29.35&	0.9694&	490 &$\alpha$=1.4& 3.96 \\
       \hline 
 Elephant 2 &$TV-L^2$  & 27.53&	20.71&	0.9009&2400&	 $\lambda$=500 & 10.64\\
   &$TV-H^{-1}$& 28.29&	21.47&	0.8997&	756&  $\lambda$=500 & 4.78   \\
      &Our&28.51&	21.70&	0.8961&	260 &$\alpha$=1.8&2.27 \\
        \hline

                Dog &$TV-L^2$  & 38.68&31.63	&0.9867 &	40265&  $\lambda$=500& 14.37\\
                       &$TV-H^{-1}$&      38.13&	31.09&	0.9758&	3946 &  $\lambda$=500 & 18.01\\
                        &Our& 39.14&32.09&	0.9764 &	830 & $\alpha$=1.6& 4.80\\
                        \hline
  
      \end{tabular}
      \caption{Comparison of Inpainting results with $TV-L^2$ and $TV-H^{-1}$ model on few Real world Images.}
 	\label{tablegray2}
 \end{table*}

 \begin{table*}[t]
 	 	\centering
 	\begin{tabular}{|c|c|c|c|c|c|c|c|} 
 		\hline
 		 Image & {\begin{tabular}[c]{@{}c@{}} Model \end{tabular}} &
 		 {\begin{tabular}[c]{@{}c@{}} PSNR \end{tabular}} & 
 		{\begin{tabular}[c]{@{}c@{}} SNR \end{tabular}}& 
 		{\begin{tabular}[c]{@{}c@{}} SSIM \end{tabular}} & 
 		{\begin{tabular}[c]{@{}c@{}} Iteration \end{tabular}}
 		&   {\begin{tabular}[c]{@{}c@{}} parameters \end{tabular}}
		&   {\begin{tabular}[c]{@{}c@{}} CPU time \end{tabular}}\\

      \hline
   Lena 1&$TV-L^2$ &31.58&	25.92&	0.9167&1091& 	 $\lambda$=100& 14.45\\
        &CVMS&31.99&	26.33&	0.9174&	645&	$\lambda$=100& 8.63\\
        &Our&          32.92&	27.27&	0.9321&	450&	$\alpha$=1.6& 9.57\\
               \hline
          Lena2  
           & $TV-L^2$ & 34.86&	29.20&	0.9456&	1990&  $\lambda$=100&26.91\\
          &CVMS&34.97&	29.31&	0.9451&	835&$\lambda$=	100& 11.30\\
            & Our& 35.62&	29.96&	0.9597&	945&	$\alpha$=1.4& 19.55\\
               \hline
            Lena3  
            & $TV-L^2$ &26.08&	20.42&	0.7972&	11550& $\lambda$	=10& 159.31\\
            &CVMS&27.67&	22.01&	0.8808&	1520	&$\lambda$=100& 21.32\\
            & Our&  28.34&	22.68&	0.8963&	1795	&$\alpha$=1.4&37.60\\         
            
                           \hline
               Barbara 1
               &$TV-L^2$  &26.93&21.04	&0.8851&	1265	& $\lambda$=100& 16.74\\
               &CVMS&27.70&	21.81&	0.8867&	970&	$\lambda$=100& 13.10\\
               & Our& 28.33&	22.45&	0.9055&	485&$\alpha$=	1.2& 10.40\\
               \hline
                Barbara 2
                &$TV-L^2$  &29.87&	23.98&	0.9346&	1890&	 $\lambda$=100& 25.72\\
                &CVMS&30.56&	24.67&	0.9340&	1020	&$\lambda$=100& 13.74\\
                 &OUR&       31.40&	25.52&	0.9505&	1055&	$\alpha$=1.2 &22.09\\
               \hline
               Barbara3 &$TV-L^2$  & 23.10&	17.21&0.6797	&7660&	 $\lambda$=10& 105.13\\
                             &CVMS&25.77&	19.88&	0.8547&	1330&	$\lambda$=100& 17.96\\
                             &OUR&          26.42&	20.53&	0.8720&	1610	&$\alpha$=1.2&35.52\\ 
               \hline
      \end{tabular}
      \caption{Comparison of  Inpainting results with $TV-L^2$ and $TV-H^{-1}$ model on  Lena and Barbara Image.}
 	\label{tablegray1}
 \end{table*}	

  \section{Convexity splitting and numerical scheme}\label{numschm}
  We will solve our proposed model \eqref{prop2} using the convexity splitting in time and Fourier spectral method in space. We replaced the Laplacian in the time discretized scheme by its fraction version and get the fractional time discretized model of the proposed model.
  Convexity splitting methods are used to solve gradient system. Split the energy functional into a convex part and a concave part and in the time discretized scheme the convex part is considered implicitly and the other one explicitly.

The proposed model \eqref{prop2} can be obtained as the gradient flow of the energies $E_1$  and  $E_2$ in $L^2$ norm where
\begin{equation}\label{f1f2}
 E_1(u)=\int_{\Omega}\Big(\frac{\mu}{2} |\Delta u|^2 + \sqrt{|\nabla u|^2+\delta^2} \Big)dx, \quad
 E_2(u)=\int_{\Omega}\frac{\lambda}{2}(f-u)^2dx.
\end{equation}
     To apply convexity splitting idea we split  $E_1$  as $E_{11}-E_{12}$ and $E_2$ as  $E_{21}-E_{22}$  where
          $$ E_{11}= \int_{\Omega} \Big( \frac{\mu}{2} |\Delta u|^2 + \frac{C_1}{2} |\nabla u|^2  \Big) dx  \quad \text{ and }$$
           $$ E_{12}= \int_{\Omega}  \Big(\frac{C_1}{2} |\nabla u|^2 -\sqrt{|\nabla u|^2+\delta^2} \Big)dx,  $$
   and
          $$  E_{21}= \int_{\Omega} \frac{C_2}{2} |u|^2 dx, \quad \text{ and } \quad E_{22}= \int_{\Omega}\Big ( -\frac{\lambda}{2}(f-u)^2 +\frac{C_2}{2} |u|^2\Big) dx,$$
            where $ C_1,C_2$ are positive and should be large enough so that $E_{ij}$ for $i,j=1,2$  becomes convex.\\
   \\ 
     So the semi-discrete time-stepping scheme for our model is
     \begin{equation}\label{eq:superposed-gradient1}
        \frac{U_{k+1}-U_k}{\Delta t} =-\nabla_{L^2}(E_{11}(U_{k+1})-E_{12}(U_k)) - \nabla_{L^2}(E_{21}(U_{k+1})- E_{22}(U_k)).
     \end{equation}
   Simplifying we get the time discretized model as:
   
   \begin{align}\label{eq:numscheme2}
      \nonumber
     \frac{U_{k+1}-U_{k}}{\Delta t} &+\mu \Delta^2 U_{k+1}-C_1\Delta U_{k+1} +C_2 U_{k+1} \\
                                                           & =\nabla \cdot \Big( \frac{\nabla U_k}{\sqrt{|\nabla U_k|^2+\delta^2}}\Big)-C_1\Delta U_k + \lambda (f-U_k) +C_2 U_k
      \end{align}
  with the boundary condition
  \begin{equation}
  \nabla U \cdot n=\nabla (\Delta U) \cdot n=0 \quad \text{on} \quad \partial \Omega.
  \end{equation}
  
        \begin{figure}
     \centering
     \begin{subfigure}{0.32\textwidth}
         \centering
         \includegraphics[width=3.9cm,height=3.5cm]{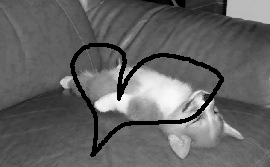}
         \caption{Damaged Dog image}
         \label{1a}
     \end{subfigure}
     \begin{subfigure}{0.32\textwidth}
         \centering
        \includegraphics[width=3.9cm,height=3.5cm]{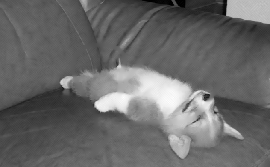}
         \caption{$\alpha=1$ (37.50)}
         \label{1b}
     \end{subfigure}
     \begin{subfigure}{0.32\textwidth}
         \centering
         \includegraphics[width=3.9cm,height=3.5cm]{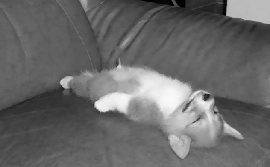}
         \caption{$\alpha=1.2$ (38.78)}
         \label{1c}
     \end{subfigure}\\
     \begin{subfigure}{0.32\textwidth}
         \centering
        \includegraphics[width=3.9cm,height=3.5cm]{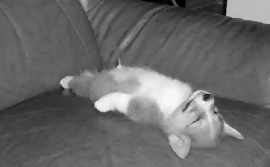}
         \caption{$\alpha=1.6$ (39.13)}
         \label{1d}
     \end{subfigure}
     \begin{subfigure}{0.32\textwidth}
         \centering
        \includegraphics[width=3.9cm,height=3.5cm]{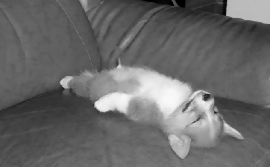}
         \caption{$\alpha=1.8$ (39.02)}
         \label{1e}
     \end{subfigure}
     \begin{subfigure}{0.32\textwidth}
         \centering
        \includegraphics[width=3.9cm,height=3.5cm]{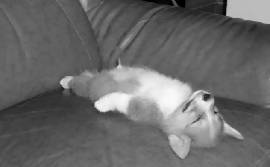}
         \caption{$\alpha$=2 (38.88)}
         \label{1f}
     \end{subfigure}
        \caption{Inpainted results of our model with different fractional power $\alpha$, tested on Dog image. In the bracket we have reported the PSNR.}
        \label{fig:dog}
  \end{figure}  
\begin{figure}
     \centering
     \begin{subfigure}{0.32\textwidth}
         \centering
         \includegraphics[width=3.9cm,height=3.5cm]{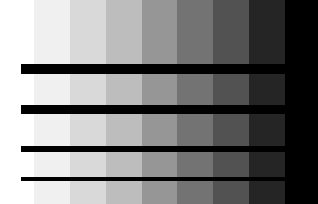}
         \caption{Damaged Grayshade}
         \label{2a}
     \end{subfigure}
     \begin{subfigure}{0.32\textwidth}
         \centering
        \includegraphics[width=3.9cm,height=3.5cm]{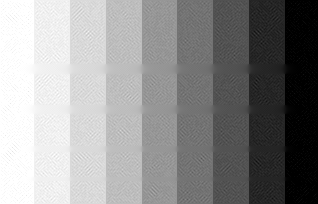}
         \caption{$\alpha=1$ (35.11)}
         \label{2b}
     \end{subfigure}
     \begin{subfigure}{0.32\textwidth}
         \centering
         \includegraphics[width=3.9cm,height=3.5cm]{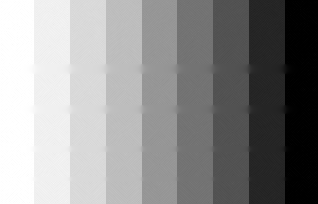}
         \caption{$\alpha=1.2$ (42.39)}
         \label{2c}
     \end{subfigure}\\
     \begin{subfigure}{0.32\textwidth}
         \centering
        \includegraphics[width=3.9cm,height=3.5cm]{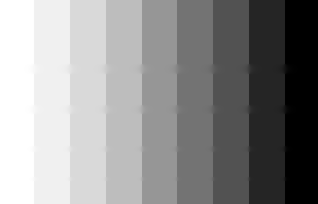}
         \caption{$\alpha=1.6$ (45.39)}
         \label{2d}
     \end{subfigure}
     \begin{subfigure}{0.32\textwidth}
         \centering
        \includegraphics[width=3.9cm,height=3.5cm]{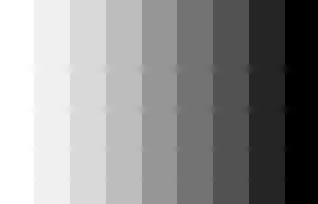}
         \caption{$\alpha=1.8$ (45.40)}
         \label{2e}
     \end{subfigure}
     \begin{subfigure}{0.32\textwidth}
         \centering
        \includegraphics[width=3.9cm,height=3.5cm]{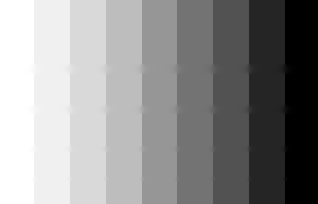}
         \caption{$\alpha$=2 (45.38)}
         \label{2f}
     \end{subfigure}
        \caption{Inpainted results of our model with different fractional power of Laplacian $\alpha$, tested on Grayshade image. In the bracket we have reported the PSNR.}
        \label{fig:gshade}
  \end{figure}
	\begin{figure}
     \centering
     \begin{subfigure}{0.35\textwidth}
         \centering
         \includegraphics[width=3.9cm,height=3.9cm]{result/inpgshade.png}
         \caption{Damaged Grayshade}
         \label{3a}
     \end{subfigure}
     \begin{subfigure}{0.35\textwidth}
         \centering
        \includegraphics[width=3.9cm,height=3.9cm]{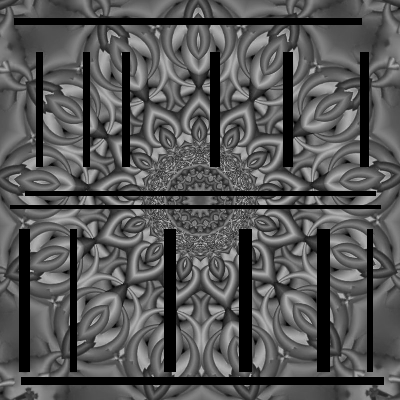}
         \caption{Damaged Kaleidoscope}
         \label{3b}
     \end{subfigure}\\
     \begin{subfigure}{0.32\textwidth}
         \centering
        \includegraphics[width=3.9cm,height=3.9cm]{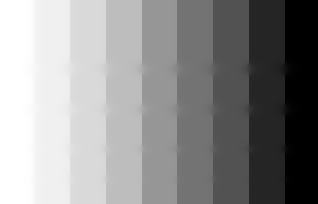}
         \caption{$TV-L^2$(41.64)}
         \label{3c}
     \end{subfigure}
     \begin{subfigure}{0.32\textwidth}
         \centering
        \includegraphics[width=3.9cm,height=3.9cm]{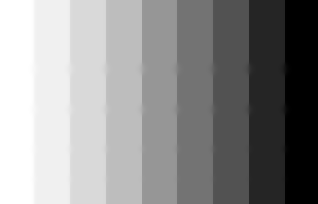}
         \caption{$TV-H^{-1}$(45.23)}
         \label{3d}
     \end{subfigure}
     \begin{subfigure}{0.32\textwidth}
         \centering
        \includegraphics[width=3.9cm,height=3.9cm]{result/fms5_gshade.png}
         \caption{Our (45.40)}
         \label{3e}
     \end{subfigure}\\
      \begin{subfigure}{0.32\textwidth}
         \centering
        \includegraphics[width=3.9cm,height=3.9cm]{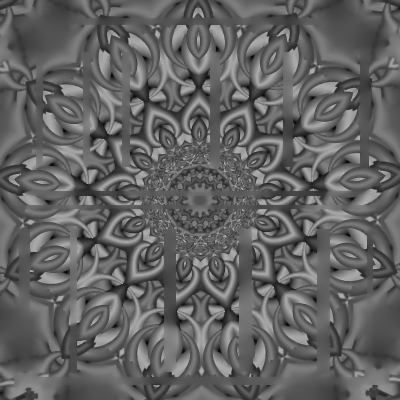}
         \caption{$TV-L^2$(26.24)}
         \label{3f}
     \end{subfigure}
      \begin{subfigure}{0.32\textwidth}
         \centering
        \includegraphics[width=3.9cm,height=3.9cm]{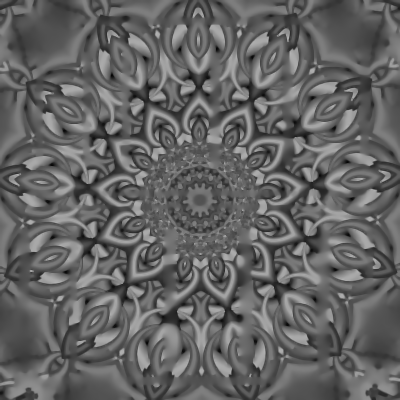}
         \caption{$TV-H^{-1}$(26.90)}
         \label{3g}
     \end{subfigure}
     \begin{subfigure}{0.32\textwidth}
         \centering
        \includegraphics[width=3.9cm,height=3.9cm]{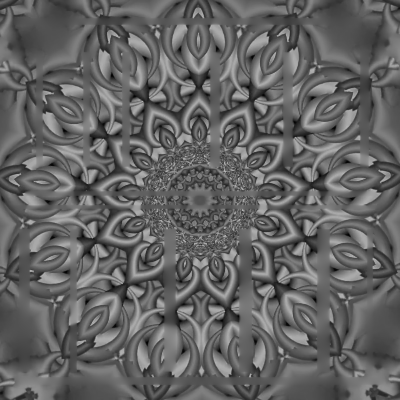}
         \caption{Our (27.89)}
         \label{3h}
     \end{subfigure}
        \caption{Comparison of our inpainted results with the results of $TV_L^2$ and $TV-H^{-1}$ model on Grayshade and Kaleidoscope image. In the bracket we have reported the PSNR.}
        \label{gshade}
\end{figure}

\begin{figure}
     \centering
     \begin{subfigure}{0.32\textwidth}
         \centering
         \includegraphics[width=3.9cm,height=3.5cm]{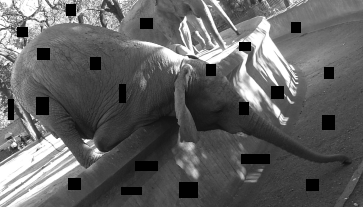}
         \caption{Damaged Elephant 1}
         \label{4a}
     \end{subfigure}
     \begin{subfigure}{0.32\textwidth}
         \centering
        \includegraphics[width=3.9cm,height=3.5cm]{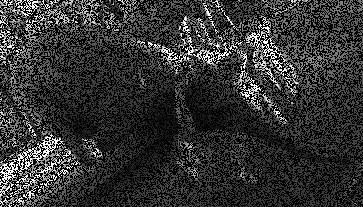}
         \caption{Damaged Elephant 2}
         \label{4b}
     \end{subfigure}
     \begin{subfigure}{0.32\textwidth}
         \centering
        \includegraphics[width=3.9cm,height=3.5cm]{result/inpdog.png}
         \caption{Damaged Dog image}
         \label{4c}
     \end{subfigure}\\
     \begin{subfigure}{0.32\textwidth}
         \centering
        \includegraphics[width=3.9cm,height=3.5cm]{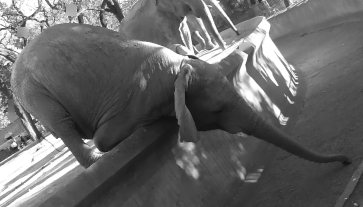}
         \caption{$TV-L^2$(35.02)}
         \label{4d}
     \end{subfigure}
     \begin{subfigure}{0.32\textwidth}
         \centering
         \includegraphics[width=3.9cm,height=3.5cm]{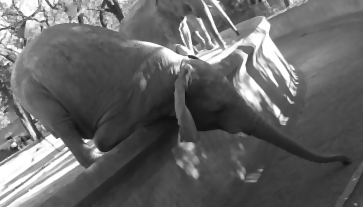}
         \caption{$TV-H^{-1}$(33.93)}
         \label{4e}
     \end{subfigure}
     \begin{subfigure}{0.32\textwidth}
         \centering
       \includegraphics[width=3.9cm,height=3.5cm]{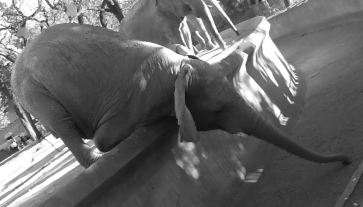}
         \caption{Our (36.16)}
         \label{4f}
     \end{subfigure}\\
      \begin{subfigure}{0.32\textwidth}
         \centering
         \includegraphics[width=3.9cm,height=3.5cm]{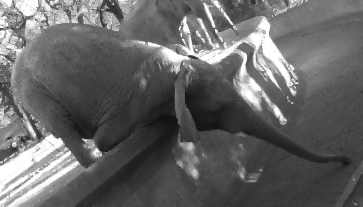}
         \caption{$TV-L^2$(27.53)}
         \label{4g}
     \end{subfigure}
      \begin{subfigure}{0.32\textwidth}
         \centering
        \includegraphics[width=3.9cm,height=3.5cm]{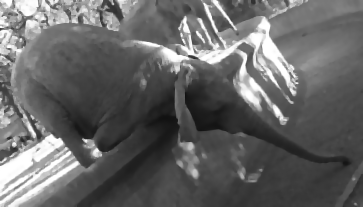}
         \caption{$TV-H^{-1}$(28.29)}
         \label{4h}
     \end{subfigure}
     \begin{subfigure}{0.32\textwidth}
         \centering
         \includegraphics[width=3.9cm,height=3.5cm]{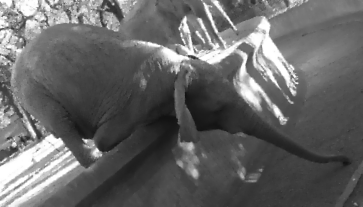}
         \caption{Our (28.51)}
         \label{4i}
     \end{subfigure}\\
      \begin{subfigure}{0.32\textwidth}
         \centering
         \includegraphics[width=3.9cm,height=3.5cm]{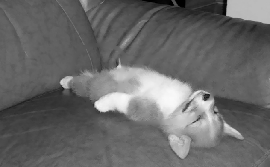}
         \caption{$TV-L^2$(38.68)}
         \label{4j}
     \end{subfigure}
      \begin{subfigure}{0.32\textwidth}
         \centering
        \includegraphics[width=3.9cm,height=3.5cm]{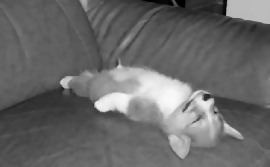}
         \caption{$TV-H^{-1}$(38.13)}
         \label{4k}
     \end{subfigure}
     \begin{subfigure}{0.32\textwidth}
         \centering
         \includegraphics[width=3.9cm,height=3.5cm]{result/fms_dog.png}
         \caption{Our (39.14)}
         \label{4l}
     \end{subfigure}
        \caption{Comparison of our inpainted results with the results of $TV_L^2$ and $TV-H^{-1}$ model on elephant and dog image. In the bracket we have reported the PSNR.}
        \label{elep}
\end{figure}
\begin{figure}
     \centering
     \begin{subfigure}{0.32\textwidth}
         \centering
        \includegraphics[width=3.9cm,height=4cm]{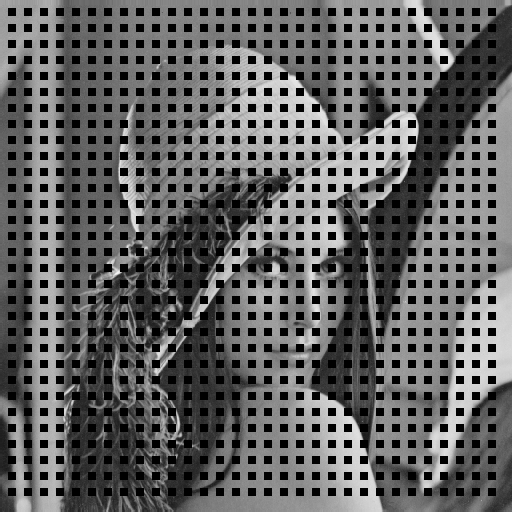}
         \caption{Damaged Lena 1}
         \label{5a}
     \end{subfigure}
     \begin{subfigure}{0.32\textwidth}
         \centering
        \includegraphics[width=3.9cm,height=4cm]{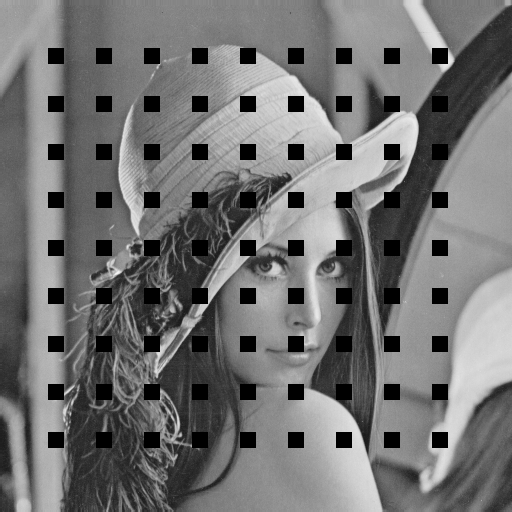}
         \caption{Damaged Lena 2}
         \label{5b}
     \end{subfigure}
     \begin{subfigure}{0.32\textwidth}
         \centering
       \includegraphics[width=3.9cm,height=4cm]{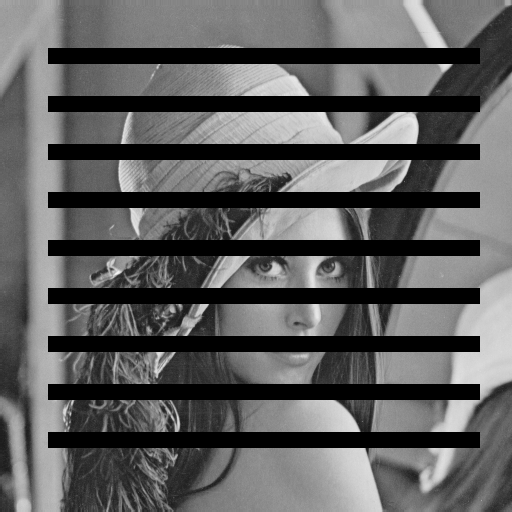}
         \caption{Damaged Lena 3}
         \label{5c}
     \end{subfigure}\\
     \begin{subfigure}{0.32\textwidth}
         \centering
        \includegraphics[width=3.9cm,height=4cm]{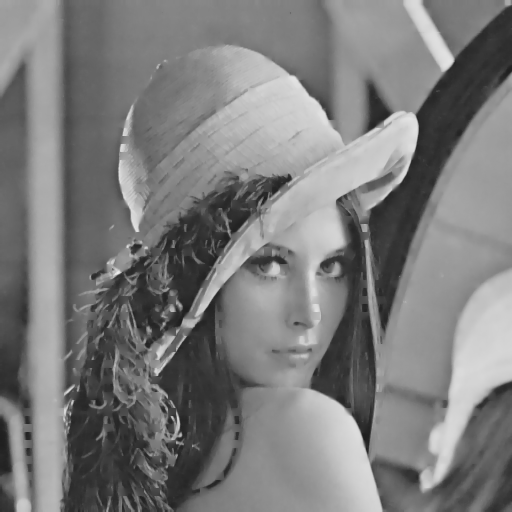}
         \caption{$TV-L^2$(31.58)}
         \label{5d}
     \end{subfigure}
     \begin{subfigure}{0.32\textwidth}
         \centering
         \includegraphics[width=3.9cm,height=4cm]{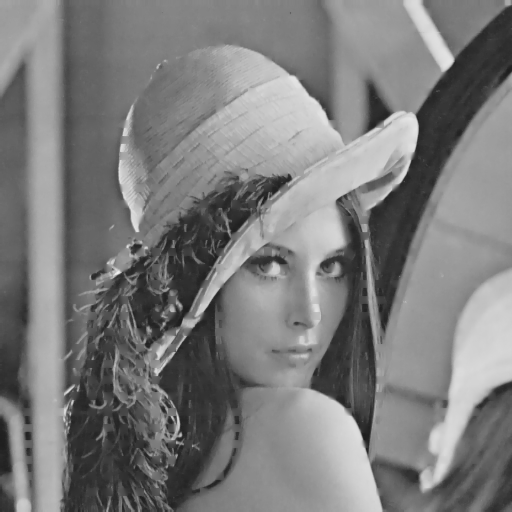}
         \caption{CVMS (31.99)}
         \label{5e}
     \end{subfigure}
     \begin{subfigure}{0.32\textwidth}
         \centering
       \includegraphics[width=3.9cm,height=4cm]{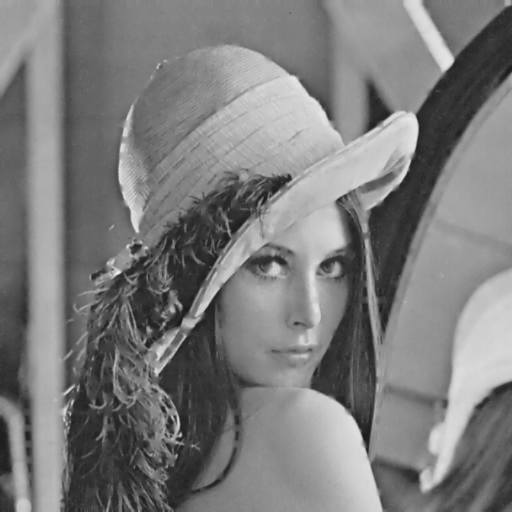}
         \caption{Our (32.92)}
         \label{5f}
     \end{subfigure}\\
      \begin{subfigure}{0.32\textwidth}
         \centering
         \includegraphics[width=3.9cm,height=4cm]{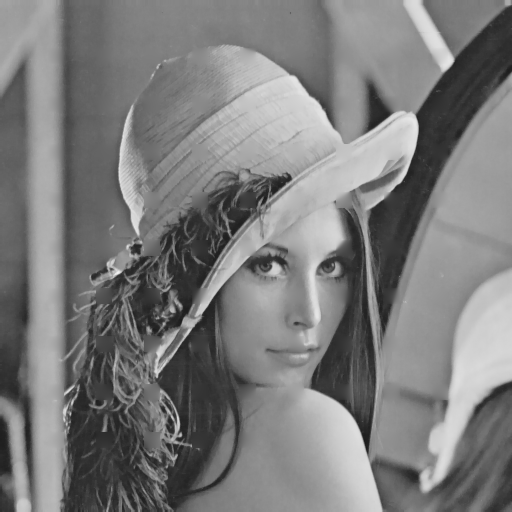}
         \caption{$TV-L^2$(34.86)}
         \label{5g}
     \end{subfigure}
      \begin{subfigure}{0.32\textwidth}
         \centering
       \includegraphics[width=3.9cm,height=4cm]{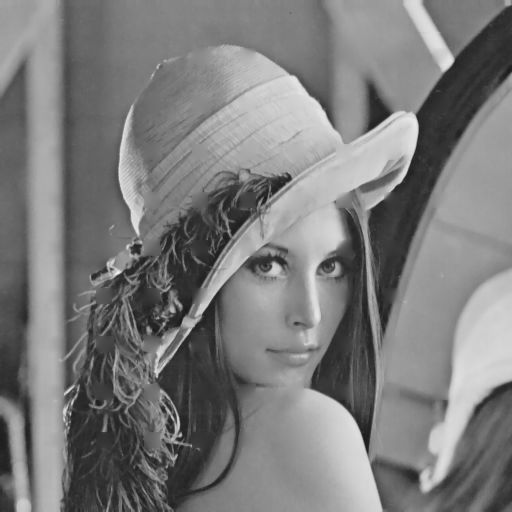}
         \caption{CVMS (34.97)}
         \label{5h}
     \end{subfigure}
     \begin{subfigure}{0.32\textwidth}
         \centering
         \includegraphics[width=3.9cm,height=4cm]{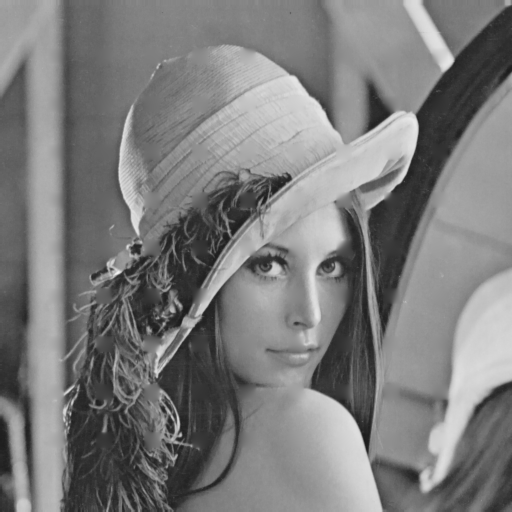}
         \caption{Our (35.62)}
         \label{5i}
     \end{subfigure}\\
      \begin{subfigure}{0.32\textwidth}
         \centering
        \includegraphics[width=3.9cm,height=4cm]{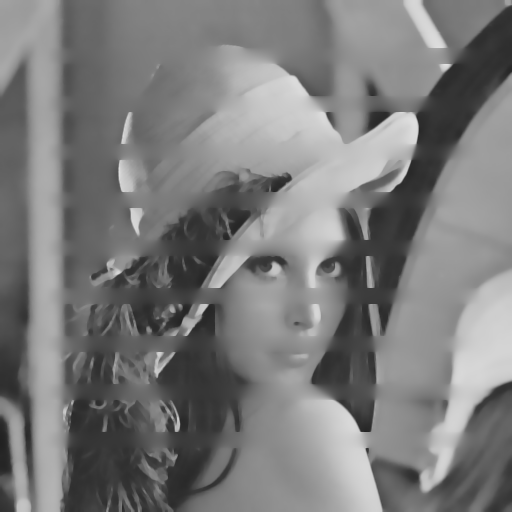}
         \caption{$TV-L^2$(26.08)}
         \label{5j}
     \end{subfigure}
      \begin{subfigure}{0.32\textwidth}
         \centering
        \includegraphics[width=3.9cm,height=4cm]{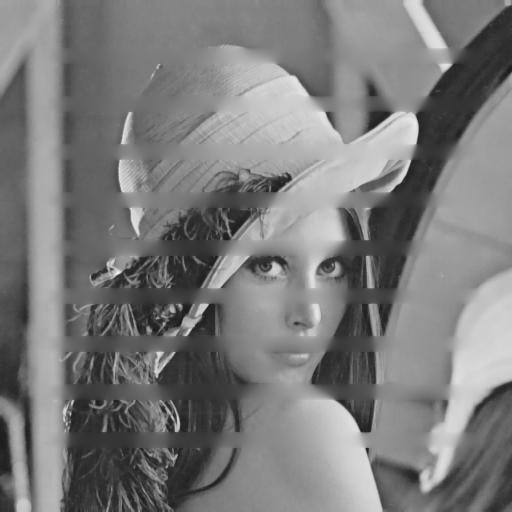}
         \caption{CVMS (27.67)}
         \label{5k}
     \end{subfigure}
     \begin{subfigure}{0.32\textwidth}
         \centering
         \includegraphics[width=3.9cm,height=4cm]{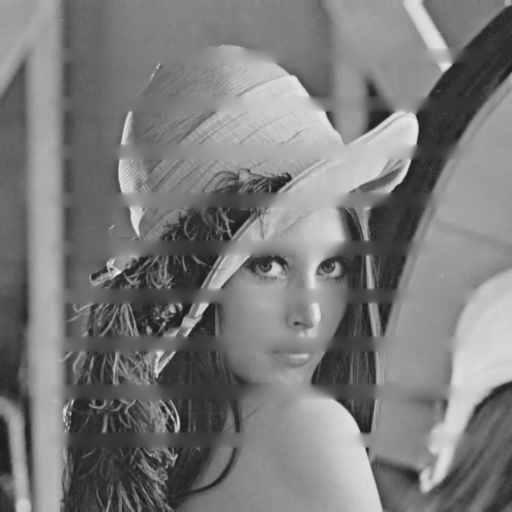}
         \caption{Our (28.34)}
         \label{5l}
     \end{subfigure}
        \caption{Comparison of our inpainted results with the results of $TV_L^2$ and $TV-H^{-1}$ model on Lena image. In the bracket we have reported the PSNR.}
        \label{lena}
\end{figure}
\begin{figure}
     \centering
     \begin{subfigure}{0.32\textwidth}
         \centering
        \includegraphics[width=3.9cm,height=4cm]{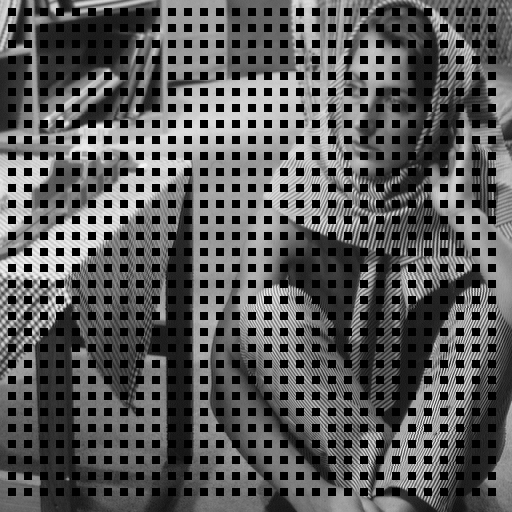}
         \caption{Damaged Barbara 1}
         \label{6a}
     \end{subfigure}
     \begin{subfigure}{0.32\textwidth}
         \centering
        \includegraphics[width=3.9cm,height=4cm]{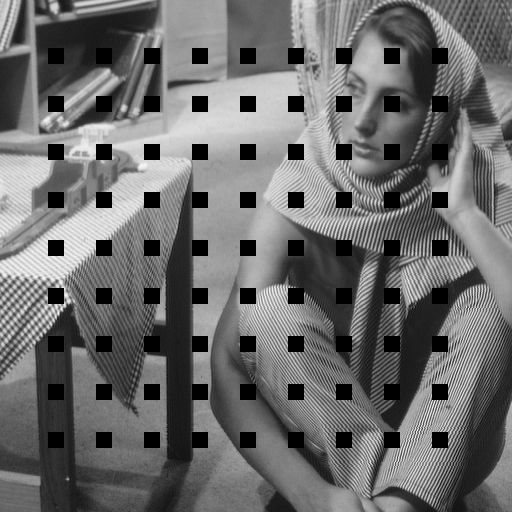}
         \caption{Damaged Barbara 2}
         \label{6b}
     \end{subfigure}
     \begin{subfigure}{0.32\textwidth}
         \centering
       \includegraphics[width=3.9cm,height=4cm]{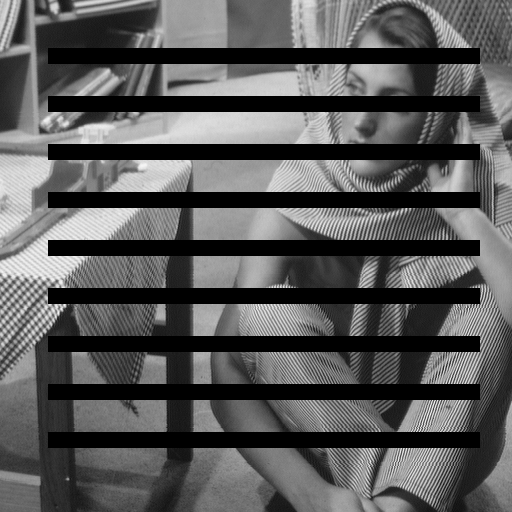}
         \caption{Damaged Barbara 3}
         \label{6c}
     \end{subfigure}\\
     \begin{subfigure}{0.32\textwidth}
         \centering
        \includegraphics[width=3.9cm,height=4cm]{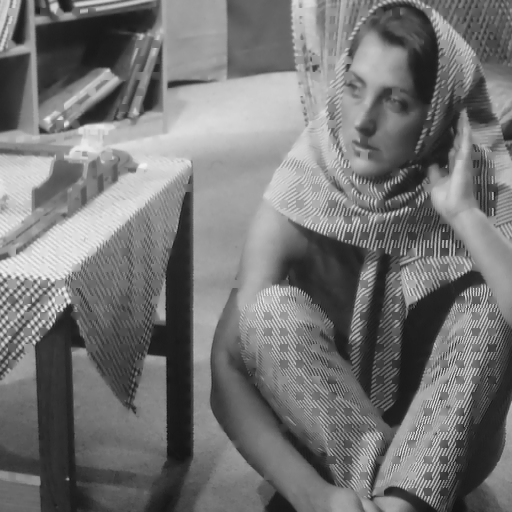}
         \caption{$TV-L^2$ (26.93)}
         \label{6d}
     \end{subfigure}
     \begin{subfigure}{0.32\textwidth}
         \centering
         \includegraphics[width=3.9cm,height=4cm]{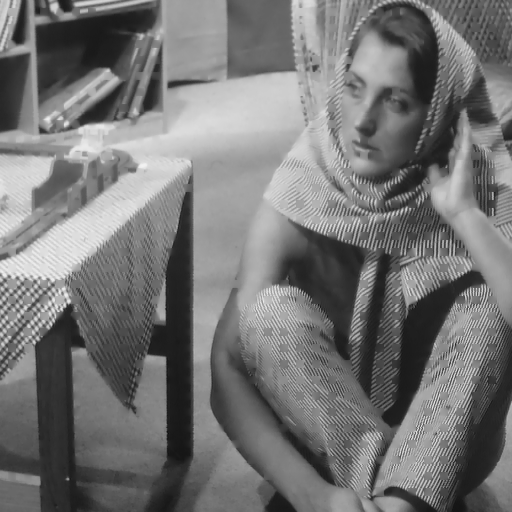}
         \caption{CVMS (27.70)}
         \label{6e}
     \end{subfigure}
     \begin{subfigure}{0.32\textwidth}
         \centering
       \includegraphics[width=3.9cm,height=4cm]{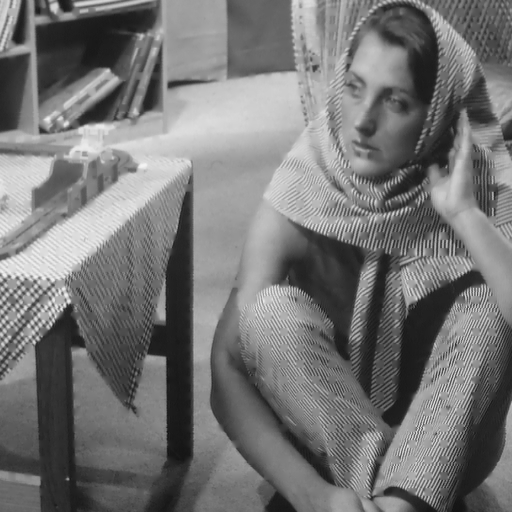}
         \caption{Our (28.33)}
         \label{6f}
     \end{subfigure}\\
      \begin{subfigure}{0.32\textwidth}
         \centering
         \includegraphics[width=3.9cm,height=4cm]{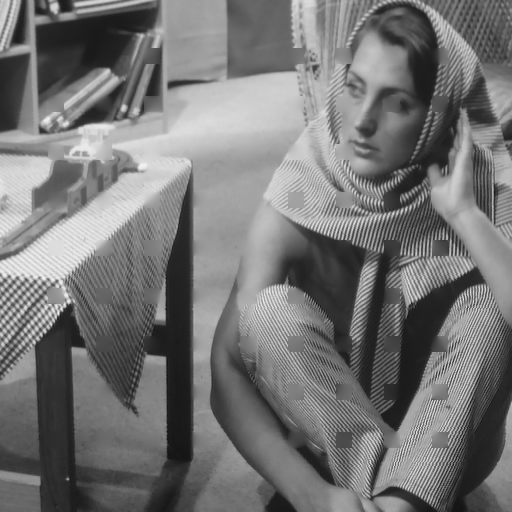}
         \caption{$TV-L^2$(29.87)}
         \label{6g}
     \end{subfigure}
      \begin{subfigure}{0.32\textwidth}
         \centering
       \includegraphics[width=3.9cm,height=4cm]{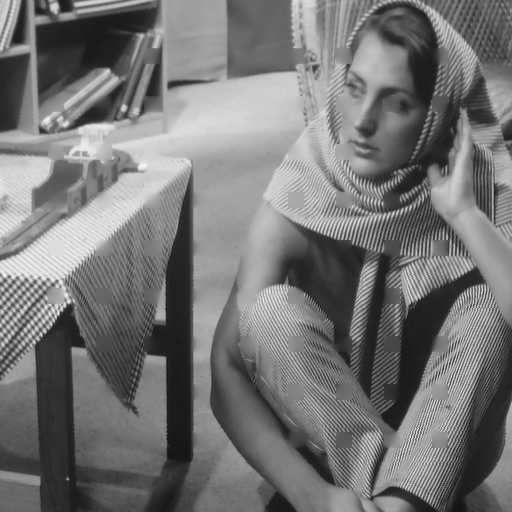}
         \caption{CVMS (30.56)}
         \label{6h}
     \end{subfigure}
     \begin{subfigure}{0.32\textwidth}
         \centering
         \includegraphics[width=3.9cm,height=4cm]{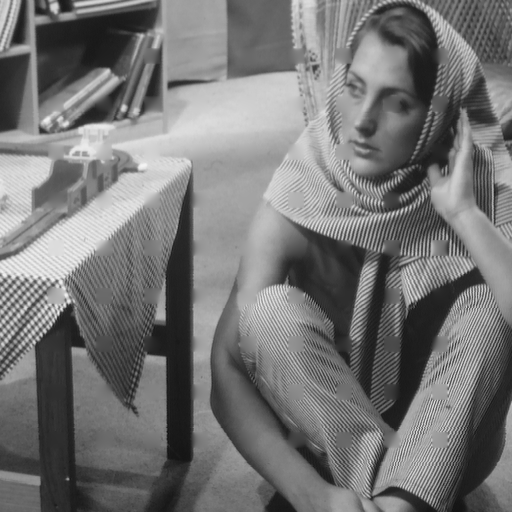}
         \caption{Our (31.40)}
         \label{6i}
     \end{subfigure}\\
      \begin{subfigure}{0.32\textwidth}
         \centering
        \includegraphics[width=3.9cm,height=4cm]{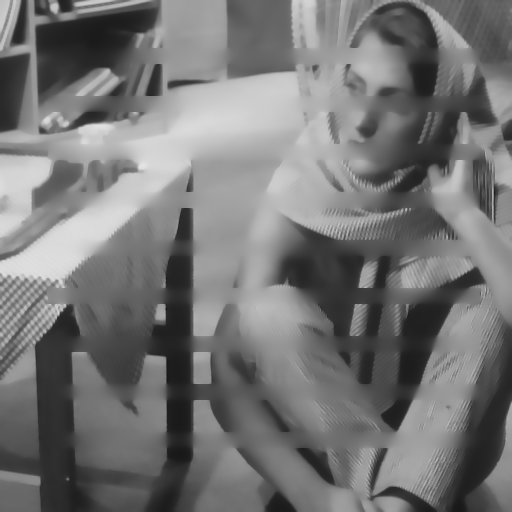}
         \caption{$TV-L^2$(23.10)}
         \label{6j}
     \end{subfigure}
      \begin{subfigure}{0.32\textwidth}
         \centering
        \includegraphics[width=3.9cm,height=4cm]{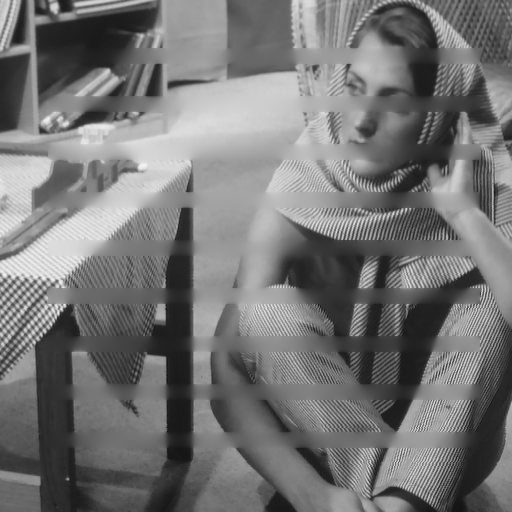}
         \caption{CVMS (25.77)}
         \label{6k}
     \end{subfigure}
     \begin{subfigure}{0.32\textwidth}
         \centering
         \includegraphics[width=3.9cm,height=4cm]{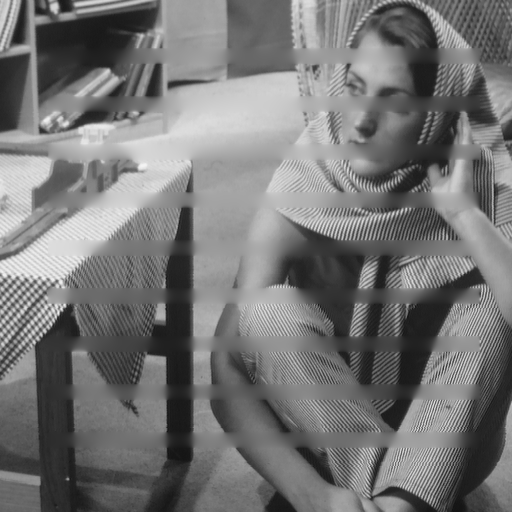}
         \caption{Our (26.42)}
         \label{6l}
     \end{subfigure}
        \caption{Comparison of inpainted results of our model with the results of $TV_L^2$ and $TV-H^{-1}$ model on Barbara image. In the bracket we have reported the PSNR.}
        \label{barb}
\end{figure}

  
Now, we replace the Laplacian $\Delta$ in \eqref{eq:numscheme2} by it's fractional counterpart $-(-\Delta)^{\frac{\alpha}{2}}$ in the spirit of~\cite{fCH} and obtained a fractional model. In general, fractional models are more effective and flexible for image processing problems~\cite{fCH}. So the equation \eqref{eq:numscheme2} reduces to:
    \begin{align}\label{eq:numscheme3}
      \nonumber
     \frac{U_{k+1}-U_{k}}{\Delta t} &+\mu (-\Delta)^{\alpha} U_{k+1}+C_1(-\Delta)^{\frac{\alpha}{2}} U_{k+1} +C_2 U_{k+1} \\
                                                           & =\nabla \cdot \Big( \frac{\nabla U_k}{\sqrt{|\nabla U_k|^2+\delta^2}}\Big)+C_1 (-\Delta)^{\frac{\alpha}{2}} U_k + \lambda (f-U_k) +C_2 U_k
      \end{align}
  with the boundary condition
  \begin{equation}
  \nabla U \cdot n=\nabla (\Delta U) \cdot n=0 \quad \text{on} \quad \partial \Omega.
  \end{equation}
 Here $\alpha$ is fractional power for $0<\alpha\leq 2.$
 We know that the Laplacian $-\Delta$ has a complete orthonormal eigenvectors $\{\varphi_{m,n}\}$ satisfying the Neumann boundary condition corresponding to eigenvalues $\lambda_{m,n}$. Consider the eigenvalue problem in $\Omega=[0,a]\times [0,b]$
 \begin{align*}
 -\Delta \varphi_{m,n}= \lambda_{m,n}\varphi_{m,n} \\
 \nabla \varphi_{m,n} \cdot n=0 \quad \text{on} \quad \partial \Omega,
 \end{align*}
 for $m,n=1,2,\dots$ where $$ \lambda_{m,n}=\pi^2 \Big( \frac{(m-1)^2}{a^2} + \frac{(n-1)^2}{b^2}\Big)$$
 $$\varphi_{m,n}=\frac{1}{\sqrt{ab}} cos\Big( \frac{(m-1)\pi x}{a}\Big)cos\Big( \frac{(n-1)\pi y}{a}\Big).$$
 
 Let us define the space \cite{fCH},
 \begin{equation}
 \mathcal{U}_{\alpha}:=\Big\{ u=\sum\limits_{m,n=1}^{\infty} \hat{u}_{m,n}\varphi_{m,n}; \hat{u}_{m,n}=\langle u,\varphi_{m,n}\rangle_{L^2}: \\ \sum\limits_{m,n=1}^{\infty} |\hat{u}_{m,n}|^2|\lambda_{m,n}|^{\frac{\alpha}{2}} <\infty, 0< \alpha\leq 2  \Big\}
 \end{equation}
 
 Then for any $u \in \mathcal{U}_{\alpha} $, the fractional Laplace operator can be defined by
 \begin{equation}\label{fracL}
 (-\Delta)^{\frac{\alpha}{2}}u=\sum\limits_{m,n=1}^{\infty} \hat{u}_{m,n}\lambda_{m,n}^{\frac{\alpha}{2}}\varphi_{m,n}.
 \end{equation}
 Hence $ (-\Delta)^{\frac{\alpha}{2}}$ has the same interpretation as $-\Delta$ in terms of its spectral decomposition. Fourier spectral methods represent the truncated series expansion when a finite number of orthonormal eigenfunction $\{\varphi_{m,n}\}$ is considered. Consistency, stabilit and convergence of the time discretized model for the case $\alpha=2$ has been carried out in the Appendix B and the analysis for the fractional one is a matter of future research as it need fractional Green's type theorem to prove the results.
Taking Fourier transform of \eqref{eq:numscheme3} and using the spectral decomposition of Laplacian \eqref{fracL} we get 
 \begin{equation}
 \widehat{U}_{k+1}(i,j)=\frac{\frac{1}{\Delta t}+(C_1\lambda_{i,j}^{\frac{\alpha}{2}} +C_2) \hat{U}_k(i,j)+\widehat{\lambda (f-U_k)}(i,j)+ \hat{\kappa}(i,j)} {\frac{1}{\Delta t}+\mu \lambda_{i,j}^{\alpha}+C_1\lambda_{i,j}^{\frac{\alpha}{2}}+C_2}
 \end{equation}
 where $\kappa =\nabla \cdot \Big ( \frac{\nabla U_k}{\sqrt{|\nabla U_k|^2+\delta^2}} \Big)$.
In every time step we will calculate the $\widehat{U_k}$ and the real part of the inverse Fourier transform will be the solution $U_k$.
 
 \section{Numerical results and discussion} \label{numresult}
 In this section, we present some numerical results on few standard test images and compared them with the
 results of  $TV-H^{-1}$ and $TV-L^2$ model. For our model, we set the parameters as $\Delta t=1$ and $\delta=0.01, \mu=0.9$, the fidelity parameter $\lambda_0=250$ and the constant occurred in convexity splitting is set as $C_1=\frac{1}{\delta}, C_2=50.$ The fractional parameter $\alpha$ is different for different images which is mentioned in the corresponding table. For comparison we have used the code of $TV-L^2$ and $TV-H^{-1}$ provided in \cite{codeinp}. To compare the quality of the results we calculate three quality metrics namely peak signal to noise ratio (PSNR), siganl to noise ratio (SNR) and structural similarity index measure (SSIM)~\cite{ssim}.  PSNR is the ratio between the maximum possible value of the original image and the mean squared error between the original and the resulting image in log scale and SNR is the ratio between the variance of the original image and the variance of the noise. SSIM represents the geometric similarity of the resulting image with the original image. For all the three quality metrics, higher the value implies better the result.
 
 To investigate the performance with different fractional power $\alpha$, we have tested our model with different values of $\alpha$ on few images.  
 In Fig.~\ref{fig:dog}, we have shown the inpainting results of dog image withe different value of $\alpha$ namely $1,1.2,1.6,1.8,2$. The results corresponding to $\alpha=1,2$ are not good enough. We have calculated the quality metrics PSNR, SNR, and SSIM for all the results and reported them in Table~\ref{table3}. From the table one can see that for $\alpha=1.4$ and $1.6$ the results are better than the others specially better than the results of the cases when $\alpha$ is an integer. Thus the model with fractional power is better than the the usual model.
 
 Similarly, we have tested our model on grayshade image  for different values of $\alpha$. The inpainted results are shown in Fig.~\ref{fig:gshade}. To compare the results we have calculated PSNR, SNR and SSIM for all the cases and reported in Table \ref{table3}. From the table one can see that in this case $\alpha=1.8$ gives the best result in this case.
 
 Now we will compare results of our model with that of $TV-H^{-1}$ and $TV-L^2$ model. In Fig.~\ref{gshade}, we have taken two images namely Grayshade image and Kaleidoscope image  and corresponding damaged image has been shown in the first row. Inpainted results of all the three models for the grayshade image are presented in second row and results of Kaleidoscope image in third row. From the figure one can see that the results of $TV-L^2$ is not good and results of $TV-H^{-1}$ and our model is better. For both the images the fractional power $\alpha$ is chosen as 1.8. PSNR, SNR and SSIM for these images are reported in Table~\ref{tablegray2}. From the table one can see that for Grayshade image our model gives little better result. For the Kaleidoscope image our model gives much better result than the results of the other two model. For example, SSIM for $TV-H^{-1}$ model is 0.87 whereas SSIM for our model is 0.91.
 
 In Fig.~\ref{elep}, we have presented the results of two images namely Elephant and dog. We have taken two different damaged region for the elephant image and one damage region for the dog image and all the damaged images are shown in the first row of Fig.~\ref{elep}. Inpainted results are shown in row two, three and four. The first column contains results of $TV-L^2$ model, second column contains the result of $TV-H^{-1}$ and results of our model is displayed in third column. For these image the fractional parameter $\alpha$ of our model is chosen as 1.4,1.8 and 1.6 respectively. From the table one can see that for all the case our model perform better than other two models in terms of PSNR, SNR and SSIM. For the Dog and first elephant image there is significant improvement in the performance in our model.
 
 In the next two experiments, we have compared our results with the results obtained by $TV-L^2$ and CVMS~\cite{twophase} model. For the CVMS model we have discretized the equation~\eqref{vms} by convexity splitting in time and Fourier spectral method in space which is presented in Appendix A.
 
 First experiment is on Lena image of size $512\times 512$. We select three different inpainting domain for Lena image and the corresponding damaged images are shown in the first row of Fig.~\ref{lena}. The results of first damaged image have been shown in second row and the results of second damaged image have been shown in third row and results of third damaged image is shown in fourth row. From the fourth row of Fig.~\ref{lena} one can see that our model is giving better result than the other two. For all the models we have reported image quality metrics PSNR, SNR, and SSIM in Table~\ref{tablegray1}. Also the best fractional power $\alpha$ is mentioned for each image in the same table. From the Table~\ref{tablegray1}, it is evident that our model perform better than the other two models in terms of PSNR, SNR, and SSIM.
 
 Finally, we did the experiment on the Barbara image of size $512\times 512$. Three different inpainting domain are chosen and the resulting degraded images are shown in the first row of Fig.~\ref{barb}. The results of $TV-L^2$ model have been shown in first column, second column contains the result of CVMS model \cite{twophase}  and our results are in third column. In this case, also we can observe from the Table~\ref{tablegray1} our model perform much better than the other two in terms of PSNR, SNR and SSIM. Note that the Barbara image contains some structures and preserving them is a challenge and in this case also our model gives better result. 
 
 \section{Conclusion}\label{conclusion}
Here we have presented a new 4th order PDE model for grayscale image inpainting using a variant of Mumford-Shah model. Convexity splitting in time has been used for time discretization and further the time discretized model has been modified by replacing the Laplace term by its fractional version.  Fourier spectral method in space has been to get the complete discretization. 
Consistency, stability and convergence have been established for the time discretized scheme. Numerical results show the superiority of our model over $TV-L^2$  and  $TV-H^{-1}$  and convex variant of Mumford-Shah model. 
\clearpage

\section*{Appendix A}
The equation \eqref{vms} can be written as gradient descent of two functional $E_1,E_2$ in $L^2$ where
\begin{equation}\label{f1f}
 E_1(u)=\int_{\Omega}\Big(\frac{\mu}{2} |\nabla u|^2 + \sqrt{|\nabla u|^2+\delta^2} \Big)dx, \quad
 E_2(u)=\int_{\Omega}\frac{\lambda}{2}(f-u)^2dx.
\end{equation}
     To apply convexity splitting idea we split  $E_1$  as $E_{11}-E_{12}$ and $E_2$ as  $E_{21}-E_{22}$  where
          $$ E_{11}= \int_{\Omega} \Big( \frac{\mu}{2} |\nabla u|^2 + \frac{C_1}{2} |\nabla u|^2  \Big) dx  \quad \text{ and }$$
           $$ E_{12}= \int_{\Omega}  \Big(\frac{C_1}{2} |\nabla u|^2 -\sqrt{|\nabla u|^2+\delta^2} \Big)dx  $$
   and
          $$  E_{21}= \int_{\Omega} \frac{C_2}{2} |u|^2 dx, \quad \text{ and } \quad E_{22}= \int_{\Omega}\Big ( -\frac{\lambda}{2}(f-u)^2 +\frac{C_2}{2} |u|^2\Big) dx,$$
            where $ C_1,C_2$ are positive and should be large enough so that $E_{ij}$ for $i,j=1,2$  becomes convex.\\
   \\ 
     So the semi-discrete time-stepping scheme for our model is
     \begin{equation}\label{eq:superposed-gradient}
        \frac{U_{k+1}-U_k}{\Delta t} =-\nabla_{L^2}(E_{11}(U_{k+1})-E_{12}(U_k)) - \nabla_{L^2}(E_{21}(U_{k+1})- E_{22}(U_k)).
     \end{equation}
     \begin{align*}
        \frac{U_{k+1}-U_k}{\Delta t} &=-\{-\mu \Delta U_{k+1}-C_1\Delta U_{k+1} +C_1\Delta U_k-\nabla \cdot \Big( \frac{\nabla U_k}{\sqrt{|\nabla U_k|^2+\delta^2}}\Big)\} \\
        & \qquad - \{C_2U_{k+1} -C_2U_k-\lambda(f-U_k)\}
     \end{align*}
     Simplifying we get the time discretized model as:
   \begin{align}\label{eq:numscheme}
      \nonumber
     \frac{U_{k+1}-U_{k}}{\Delta t} &-\mu \Delta U_{k+1}-C_1\Delta U_{k+1} +C_2 U_{k+1} \\
                                                           & =\nabla \cdot \Big( \frac{\nabla U_k}{\sqrt{|\nabla U_k|^2+\delta^2}}\Big)-C_1\Delta U_k + \lambda (f-U_k) +C_2 U_k
      \end{align}
To get the complete discretized scheme we will use Fourier spectral method in space.  Let $\widehat{U_k}$ be the discrete Fourier trans form of $U_k$. Taking Fourier transform of \eqref{eq:numscheme} and using the relation $\widehat{\Delta u}=L\hat{u}$ and rearranging we get the final numerical scheme for the first model as:
 \begin{equation}
 \widehat{U}_{k+1}(i,j)=\frac{(\frac{1}{\Delta t}-C_1L_{i,j}+C_2)\widehat{U_k}(i,j)+\widehat{\kappa_k}(i,j)+\widehat{\lambda(f-U_k)}(i,j)}{\frac{1}{\Delta t}-(\mu +C_1)L_{i,j}+C_2}
 \end{equation}
 where $\kappa_k=\nabla\cdot\frac{\nabla U_k}{\sqrt{|\nabla U_k|^2+\delta^2}}$ and $L=(L_{i,j})$ is matrix with size $m\times n$ with 
 $$L_{i,j}=\frac{2}{(\Delta x)^2} \Big(cos (\frac{2\pi i}{m})-1\Big)+\frac{2}{(\Delta y)^2} \Big(cos (\frac{2\pi j}{n})-1\Big).$$
        
   \section*{Appendix B}     
 Now, we will discuss the consistency, stability and convergence of the time discretized model \eqref{eq:numscheme2} that is for the fractional model with $\alpha=2.$
         \begin{theorem}[Consistency, Stability and Convergence]
       Let $u$ be the exact solution of \eqref{prop2} and $u_k=u(k\Delta t)$ be the exact solution at time $k\Delta t$ for a time step $\Delta t > 0$ and $k\in \mathbb{N}.$
 Let $U_k$ be the $k$-th iterate of\eqref{eq:numscheme2} with constant $C1>\frac{1}{\delta},C_2>\lambda_0$. Then the following statements are true:
 \begin{enumerate}[(a)]
\item Under the assumption that $\|u_{tt}\|_2$ and $\|u_t\|_2$ are bounded, the numerical scheme  is consistent with the continuous
 equation  and of order one in time.
\item the solution $U_k$ is bounded on a finite time interval $[0,T]$, for all $\Delta t>0$. In particular for $k\Delta t\leq T$, $T>0$ fixed, we have for every $\Delta t>0$
 	\begin{align}
 		& \|\nabla U_k\|^2+\Delta t K_1\|\Delta U_{k+1}\|^2+\Delta t K_2 \|\nabla\Delta U_k\|^2 \\ \nonumber
         & \leq e^{K T}\Big(\|\nabla U_0\|^2 +\Delta t K_1\|\Delta U_{0}\|^2+\Delta t K_2\|\nabla\Delta U_0\|^2+ \Delta t TC(\Omega,D,\lambda_0,f)\Big),
 	\end{align}
	 \item The discretization error $e_k=u_k-U_k$. For smooth solution $u_k$ and $U_k$, the error $e_k$ converges to $0$ as $\Delta t\to 0$.
 In particular we have for  $k\Delta t\leq T$, $T>0$ fixed, that
 	\begin{center}
 		$\|\nabla e_k\|^2+\Delta t M_1\|\Delta U_{k+1}\|^2+\Delta t M_2\|\nabla\Delta e_k\|^2 \leq \frac{T}{M_3}e^{M_4T}(\Delta t)^2$
 	\end{center}
	for suitable constants $M1,M_2,M_3,M4.$
 \end{enumerate}	
      \end{theorem}
    \proof (a)
  Putting $u=u_k$ in equation \eqref{prop2} and $U_k=u_k$ in \eqref{eq:numscheme2} and subtracting the later we get the error $\tau_k$ at $k^{th}$ time step as
    \begin{equation}\label{lte}
    \tau_k=\frac{u_{k+1}-u_k}{\Delta t}-u_t(k\Delta t)+\mu \Delta^2 (u_{k+1}-u_k)-C_1 \Delta (u_{k+1}-u_k)+C_2(u_{k+1}-u_k)
    \end{equation}
     Then applying Taylors theorem we get,
    $$\tau_k=(u_t(k\Delta t)+\Delta t u_{tt}(\xi_1))-u_t(k\Delta t)+ \Delta t (\mu  \Delta^2 u_t(\xi_2)-C_1  \Delta u_t(\xi_3)+C_2 u_t(\xi_4))$$
    Now using the boundedness assumptions we get $\|\tau_k\|=O(\Delta t)$. Hence the global truncation error $\tau=\max\limits_k \|\tau_k\|=O(\Delta t)$.
 \endproof
 
    \proof  (b)
   Consider the time discretized scheme
   \begin{align}\label{timeschm}
      \nonumber
     \frac{U_{k+1}-U_{k}}{\Delta t} &+\mu \Delta^2 U_{k+1}-C_1\Delta U_{k+1} +C_2 U_{k+1} \\
                                                           & =\nabla \cdot \Big( \frac{\nabla U_k}{\sqrt{|\nabla U_k|^2+\delta^2}}\Big)-C_1\Delta U_k + \lambda (f-U_k) +C_2 U_k
      \end{align}
      Multiplying the time discretised scheme \eqref{timeschm} by $-\Delta U_{k+1}$ and integrating over $\Omega$ we get
      \begin{align*}
      \frac{1}{\Delta t} \Big( \|\nabla U_{k+1}\|^2-&\langle \nabla U_k, U_{k+1}\rangle \Big)+ \mu \|\nabla \Delta U_{k+1}\|^2+C_1\|\Delta U_{k+1}\|^2 +C_2\|\nabla U_{k+1}\|^2\\
      &=\Big\langle \nabla\nabla \cdot \frac{\nabla U_k}{\sqrt{|\nabla U_k|^2+\delta^2}}, \nabla U_{k+1}  \Big\rangle +C_1 \langle \Delta U_k,\Delta U_{k+1} \rangle\\
     & \qquad+\langle \nabla \lambda(f-U_k),\nabla U_{k+1} \rangle +C_2  \langle \nabla U_k,\nabla U_{k+1} \rangle
       \end{align*}
       
    Using Young's inequality on the product terms  we get,
    \begin{align}\label{afteryng}
    \nonumber
    &  \frac{1}{2\Delta t} \Big( \|\nabla U_{k+1}\|^2-\|\nabla U_k\|^2 \Big)+ \mu \|\nabla \Delta U_{k+1}\|^2+C_1\|\Delta U_{k+1}\|^2 +C_2\|\nabla U_{k+1}\|^2\\ \nonumber
      &\quad \leq \delta_1 \| \nabla U_{k+1}\|^2+\frac{1}{4\delta_1}\|  \nabla\nabla \cdot \frac{\nabla U_k}{\sqrt{|\nabla U_k|^2+\delta^2}}\|^2
      +C_1\delta_2 \| \Delta U_{k+1} \|^2+\frac{C_1}{4\delta_2}\| \Delta U_k\|^2\\ 
      & \quad +\delta_3 \| \nabla U_{k+1} \|^2+\frac{1}{4\delta_3}\| \lambda(f-U_k) \|^2+C_2\delta_4 \| \nabla U_{k+1} \|^2+\frac{C_2}{4\delta_4}\| \nabla U_k\|^2
    \end{align}
     Using the following two estimates \cite{nondenoising} in \eqref{afteryng}
     \begin{equation}\label{bnd1}
     \| \lambda(f-U_k) \|^2 \leq 2 \lambda_0^2 \|\nabla U_k\|^2+C(\Omega,D,\lambda_0,f)
     \end{equation}
      and 
        \begin{equation}\label{bnd2}
        \|  \nabla\nabla \cdot \frac{\nabla U_k}{\sqrt{|\nabla U_k|^2+\delta^2}}\|^2 \leq \bar{C}(\delta,\Omega)(\| \nabla U_k\|^2+\|\nabla \Delta U_k\|^2)
        \end{equation}
        we get,
        \begin{align*}
       & \Big(\frac{1}{2\Delta t} +C_2(1-\delta_4)-(\delta_1+\delta_3)\Big) \|\nabla U_{k+1}\|^2+C_1(1-\delta_2) \|\Delta U_{k+1}\|^2+ \mu \|\nabla \Delta U_{k+1}\|^2\\
        & \quad \leq  \Big(\frac{1}{2\Delta t} +\frac{\bar{C}}{4\delta_1}+\frac{C_2}{4\delta_4}+\frac{2\lambda_0^2}{4\delta_3}\Big) \|\nabla U_{k}\|^2+\frac{C_1}{4\delta_2}\| \Delta U_k\|^2 +\frac{\bar{C}}{4\delta_1}\|\nabla \Delta U_{k+1}\|^2+C(\Omega,D,\lambda_0,f)
        \end{align*}
       Let us choose $\delta_1=\delta_3=\frac{1}{4}$ and $\delta_2=\delta_4=\frac{1}{2}$, we get,
       \begin{align*}
        & \Big(\frac{1}{2\Delta t} +\frac{C_1-1}{2}\Big) \|\nabla U_{k+1}\|^2+\frac{C_1}{2} \|\Delta U_{k+1}\|^2+ \mu \|\nabla \Delta U_{k+1}\|^2\\
        & \quad \leq  \Big(\frac{1}{2\Delta t} +\bar{C}+\frac{C_2}{2}+2\lambda_0^2\Big) \|\nabla U_{k}\|^2+\frac{C_1}{2}\| \Delta U_k\|^2 +\bar{C}\|\nabla \Delta U_{k+1}\|^2+C(\Omega,D,\lambda_0,f)
       \end{align*}
       Since $C_2>1$, all the coefficients are positive. Now multiplying both sides by $2 \Delta t$ we get,
       \begin{align}\nonumber
     &  A_1  \|\nabla U_{k+1}\|^2 +\Delta t C_1 \|\Delta U_{k+1}\|^2+ 2\Delta t \mu \|\nabla \Delta U_{k+1}\|^2\\
     & \quad \leq A_2  \|\nabla U_{k}\|^2 +\Delta t C_1 \|\Delta U_{k}\|^2+ 2\Delta t \bar{C} \|\nabla \Delta U_{k}\|^2+2 \Delta t C(\Omega,D,\lambda_0,f), 
       \end{align}
       where $A_1=1+\Delta t(C_2-1), A_2=1+ \Delta t(C_2+2\bar{C}+4 \lambda_0^2)$.
       Dividing both sides by $A_1$ we get,
        \begin{align}\nonumber
     &   \|\nabla U_{k+1}\|^2 +\Delta t \frac{C_1}{A_1} \|\Delta U_{k+1}\|^2+ \Delta t \frac{2\mu}{A_1} \|\nabla \Delta U_{k+1}\|^2\\
     & \quad \leq \frac{A_2}{A_1} \Big( \|\nabla U_{k}\|^2 +\Delta t \frac{C_1}{A_2} \|\Delta U_{k}\|^2+ \Delta t \frac{2\bar{C}}{A_2} \|\nabla \Delta U_{k}\|^2\Big)+\frac{2 \Delta t }{A_1} C(\Omega,D,\lambda_0,f), 
       \end{align}
Choose $A_1,A_2$ in such a way that $\frac{A_2}{A_1}> 1$ and $\frac{\mu A_2}{\bar{C}A_1}> 1$ (this is possible if we choose big enough $\lambda_0$), then above equation reduces to
\begin{align}\nonumber
     &   \|\nabla U_{k+1}\|^2 +\Delta t \frac{C_1}{A_1} \|\Delta U_{k+1}\|^2+ \Delta t \frac{2\mu}{A_1} \|\nabla \Delta U_{k+1}\|^2\\
     & \quad \leq \frac{A_2}{A_1} \Big( \|\nabla U_{k}\|^2 +\Delta t \frac{C_1}{A_1} \|\Delta U_{k}\|^2+ \Delta t \frac{2\mu}{A_1} \|\nabla \Delta U_{k}\|^2\Big)+ \Delta t C(\Omega,D,\lambda_0,f).
     \end{align}
     So by induction we will get,
     \begin{align}\nonumber
     &   \|\nabla U_{k+1}\|^2 +\Delta t \frac{C_1}{A_1} \|\Delta U_{k+1}\|^2+ \Delta t \frac{2\mu}{A_1} \|\nabla \Delta U_{k+1}\|^2\\
     & \quad \leq \Big(\frac{A_2}{A_1}\Big)^k \Big( \|\nabla U_{0}\|^2 +\Delta t \frac{C_1}{A_1} \|\Delta U_{0}\|^2+ \Delta t \frac{2\mu}{A_1} \|\nabla \Delta U_{0}\|^2\Big)+ \Delta t \sum\limits_{i=0}^{k-1}\Big(\frac{A_2}{A_1}\Big)^iC(\Omega,D,\lambda_0,f).
     \end{align}
     Expressing $\frac{A_2}{A_1}$ as =$1+K_3 \Delta t$ and using the inequality $(1+nx)^k\leq e^{kx}$, we get for $k\Delta t< T$,
     \begin{align}\nonumber
     & \|\nabla U_{k+1}\|^2 +\Delta t \frac{C_1}{A_1} \|\Delta U_{k+1}\|^2+ \Delta t \frac{2\mu}{A_1} \|\nabla \Delta U_{k+1}\|^2\\
     & \quad \leq e^{K_3T} \Big( \|\nabla U_{0}\|^2 +\Delta t \frac{C_1}{A_1} \|\Delta U_{0}\|^2+ \Delta t \frac{2\mu}{A_1} \|\nabla \Delta U_{0}\|^2+TC(\Omega,D,\lambda_0,f)\Big).
     \end{align}
   \endproof

 \begin{proof}(c)
 Putting the value of $u_t(k\Delta t)$ from \eqref{prop2} in equation \eqref{lte} and simplifying we will have,
 \begin{align}\label{Lte1}
 \nonumber
 \tau_k=\frac{u_{k+1}-u_k}{\Delta t}&-\nabla \cdot \frac {\nabla u_k}{\sqrt{|\nabla u_k|^2+\delta^2}} -\lambda(f-u_k)+\mu \Delta^2 u_{k+1}\\ 
 &  -C_1 \Delta (u_{k+1}-u_k)+C_2(u_{k+1}-u_k)
 \end{align}
 Now subtracting \eqref{timeschm} from \eqref{Lte1} and putting $e_k=u_k-U_k$ we get,
 \begin{align}\label{eq:err}
 \nonumber
& \frac{e_{k+1}-e_k}{\Delta t}+\mu \Delta^2e_{k+1}-C_1\Delta e_{k+1}+C_2e_{k+1}\\
 &=\Big( \nabla \cdot \frac {\nabla u_k}{\sqrt{|\nabla u_k|^2+\delta^2}}-\nabla \cdot \frac {\nabla U_k}{\sqrt{|\nabla U_k|^2+\delta^2}}\Big)-C_1\Delta e_k+C_2e_k-\lambda e_k +\tau_k
 \end{align}
 Multiplying equation \eqref{eq:err} by $e_{k+1}$ and integrating over $\Omega$ and arguing as proof (b) and following \cite{tvh} we will get the
required result.
 \end{proof}

\end{document}